# SPATIAL AGGREGATION OF LOCAL LIKELIHOOD ESTIMATES WITH APPLICATIONS TO CLASSIFICATION[1]

By Denis Belomestny and Vladimir Spokoiny

*Weierstrass Institute*

This paper presents a new method for spatially adaptive local (constant) likelihood estimation which applies to a broad class of nonparametric models, including the Gaussian, Poisson and binary response models. The main idea of the method is, given a sequence of local likelihood estimates ("weak" estimates), to construct a new aggregated estimate whose pointwise risk is of order of the smallest risk among all "weak" estimates. We also propose a new approach toward selecting the parameters of the procedure by providing the prescribed behavior of the resulting estimate in the simple parametric situation. We establish a number of important theoretical results concerning the optimality of the aggregated estimate. In particular, our "oracle" result claims that its risk is, up to some logarithmic multiplier, equal to the smallest risk for the given family of estimates. The performance of the procedure is illustrated by application to the classification problem. A numerical study demonstrates its reasonable performance in simulated and real-life examples.

**1. Introduction.** This paper presents a new method of spatially adaptive nonparametric estimation based on the aggregation of a family of local likelihood estimates. As a main application of the method, we consider the problem of building a classifier on the base of the given family of $k$-NN or kernel classifiers.

The local likelihood approach has been intensively discussed in recent years; see, for example, Tibshirani and Hastie [16], Staniswalis [15] and Loader [11]. We refer to Fan, Farmen and Gijbels [5] for a nice and detailed overview of the local maximum likelihood approach and related lit-

Received June 2005; revised January 2007.

[1]Supported by the Deutsche Forschungsgemeinschaft through the SFB 649 Economic Risk.

erature. Similarly to nonparametric smoothing in the regression or density framework, an important issue for local likelihood modeling is the choice of localization (smoothing) parameters. Different types of model selection techniques based on the asymptotic expansion of the local likelihood are mentioned in Fan, Farmen and Gijbels [5], which include global, as well as variable, bandwidth selection. However, the finite sample performance of estimators based on bandwidth or model selection is often rather unstable; see, for example, Breiman [2]. This point is particulary critical for the local or pointwise model selection procedures like Lepski's method. In spite of the nice theoretical properties (see Lepski, Mammen and Spokoiny [8], Lepski and Spokoiny [9] or Spokoiny [14]), the resulting estimates suffer from high variability due to a pointwise model choice, especially for a large noise level. This suggests that in some cases, the attempt to identify the true model is not necessarily appropriate. One approach to reducing variability in adaptive estimation is model mixing or aggregation. Catoni [4] and Yang [18], among others, have suggested global aggregating procedures that achieve the minimal estimation risks over the family of given "weak" estimates. In the regression setting, Juditsky and Nemirovski [7] have developed aggregation procedures which have a risk within a multiple of the smallest risk in the class of all convex combinations of "weak" estimates plus $\log(n)/n$. Tsybakov [17] has discussed asymptotic minimax rates for aggregation. Aggregation for density estimation has been studied by Li and Barron [10] and more recently by Rigollet and Tsybakov [13]. To the best of our knowledge, pointwise aggregation has not yet been considered.

Our approach is based on the idea of the spatial (pointwise) aggregation of a family of local likelihood estimates ("weak" estimates) $\widetilde{\theta}^{(k)}$. The main idea is, given the sequence $\{\widetilde{\theta}^{(k)}\}$, to construct in a data-driven way, for every point $x$, the "optimal" aggregated estimate $\widehat{\theta}(x)$. "Optimality" means that this estimate satisfies some kind of "oracle" inequality, that is, its pointwise risk does not exceed the smallest pointwise risk among all "weak" estimates up to a logarithmic multiple.

Our algorithm can be roughly described as follows. Let $\{\widetilde{\theta}^{(k)}(x)\}$, $k = 1, \ldots, K$, be a sequence of "weak" local likelihood estimates at a point $x$, ordered according to their variability, which decreases with $k$. Starting with $\widehat{\theta}^{(1)}(x) = \widetilde{\theta}^{(1)}(x)$, an aggregated estimate $\widehat{\theta}^{(k)}(x)$ at any step $1 < k \leq K$ is constructed by mixing the previously constructed aggregated estimate $\widehat{\theta}^{(k-1)}(x)$ with the current "weak" estimate $\widetilde{\theta}^{(k)}(x)$,

$$\widehat{\theta}^{(k)}(x) = \gamma_k \widetilde{\theta}^{(k)}(x) + (1 - \gamma_k)\widehat{\theta}^{(k-1)}(x),$$

and taking $\widehat{\theta}^{(K)}(x)$ as a final estimate. The mixing parameter $\gamma_k$ (which may depend on the point $x$) is defined using a measure of statistical difference between $\widehat{\theta}^{(k-1)}(x)$ and $\widetilde{\theta}^{(k)}(x)$. In particular, $\gamma_k$ is equal to zero if $\widehat{\theta}^{(k-1)}(x)$ lies outside the confidence set around $\widetilde{\theta}^{(k)}(x)$. In view of the sequential and pointwise nature of the algorithm, the suggested procedure is called *Spatial Stagewise Aggregation* (SSA). Important features of the proposed procedure are its simplicity and applicability to a variety of problems including Gaussian, binary, Poisson regression, density estimation, classification, etc. The procedure does not require any splitting of the sample as many other aggregation procedures do (cf. Yang [18]). Furthermore, the theoretical properties of SSA can be rigorously studied. In particular, we establish precise nonasymptotic "oracle" results which are applicable under very mild conditions in a rather general setting. We also show that the oracle property automatically implies spatial adaptivity of the proposed estimate.

Another important feature of the procedure is that it can be easily implemented and the problem of selecting the tuning parameters can be carefully addressed.

Our simulation study confirms nice finite sample performance of the procedure for a broad class of models and problems. We only show the results for the classification problem as the most interesting and difficult one. Some more examples for univariate regression and density estimation can be found in our preprint Belomestny and Spokoiny [1]. Section 4 shows how the SSA procedure can be applied to aggregating kernel and $k$-NN classifiers in the classification problem. Although these two nonparametric classifiers are rather popular, the problem of selecting the smoothing parameter (the bandwidth for the kernel classifier or the number of neighbors for the $k$-NN method) has not yet been satisfactorily addressed. Again, the SSA-based classifier demonstrates the "oracle" quality in terms of the both pointwise and global misclassification errors. This application clearly shows one more important feature of the SSA method: it can be applied to an arbitrary design with design space of arbitrary dimension. This is illustrated by simulated and real-life classification examples in dimensions up to 10.

The procedure proposed in this paper is limited to aggregating the kernel-type estimates which are based on local constant approximation. The modern statistical literature usually considers the more general local linear (polynomial) approximation of the underlying function. However, for this paper, we have decided for several reasons to restrict our attention to the local constant case. The most important one is that for the examples and applications we consider in this paper, the use of the local linear methods does not improve (and even degrades) the quality of estimation. Our experience strongly confirms that for problems like classification, local constant smoothing combined with the aggregation technique delivers reasonable finite sample quality.

Our theoretical study is split into two major parts. Section 2 introduces the local parametric setting to be considered and extends the parametric risk bounds to the local parametric and nonparametric situation under the so-called "small modeling bias" condition. The main result (Corollary 2.6) claims that the parametric risk bounds continue to apply provided that this condition is fulfilled. One possible interpretation of our adaptive procedure is the search for the largest localizing scheme for which the "small modeling bias" condition still holds. Theoretical properties of the aggregation procedure are presented in Section 5. The main result states the "oracle" property of the SSA estimate: the risk of the aggregated estimate is, within a log-multiple, as small as the risk of the best "weak" estimate for the function under consideration. The results are established in the precise nonasymptotic way for a rather general likelihood setting under mild regularity conditions. Moreover, our approach establishes a link between parametric and nonparametric theory. In particular, we show that the proposed method delivers root-$n$ accuracy in the parametric situation. In the nonparametric case, the quality corresponds to the best parametric approximation. Both the theoretical study and the motivation of the procedure employ some exponential bounds for the likelihood which are given in Section 2.2. An important feature of our theoretical study is that the problem of selecting the tuning parameters is also discussed in detail. We offer a new approach, in which the parameters of the procedure are selected to provide the desirable performance of the method in the simple parametric situation. This is similar to the hypothesis problem approach when the critical values are selected using the performance of the test statistic under the simple null hypothesis; see Section 3.3.1 for a detailed explanation.

**2. Local constant modeling for exponential families.** This section presents some results on local constant likelihood estimation. We begin by describing the model under consideration. Suppose we are given independent random data $Z_1, \ldots, Z_n$ of the form $Z_i = (X_i, Y_i)$. Here, every $X_i$ denotes a vector of "features" or explanatory variables which determine the distribution of the "observation" $Y_i$. For simplicity, we assume that the $X_i$'s are valued in the finite-dimensional Euclidean space $\mathscr{X} = \mathbb{R}^d$ and that the $Y_i$'s belong to $\mathbb{R}$. The vector $X_i$ can be viewed as a location and $Y_i$ as the "observation at $X_i$." Our model assumes that the distribution of each $Y_i$ is determined by a finite-dimensional parameter $\theta$ which may depend on the location $X_i$.

More precisely, let $\mathscr{P} = (P_\theta, \theta \in \Theta)$ be a parametric family of distributions dominated by a measure $P$. In this paper, we only consider the case when $\Theta$ is a subset of the real line. By $p(\cdot, \theta)$ we denote the corresponding density. We consider the regression-like model in which every "response" $Y_i$ is, conditionally on $X_i = x$, distributed with the density $p(\cdot, f(x))$ for some unknown function $f(x)$ on $\mathscr{X}$ with values in $\Theta$. The model under consideration can be written as

$$Y_i \sim P_{f(X_i)}.$$

The aim of the data analysis is to estimate the function $f(x)$. For related models, see Fan and Zhang [6] and Cai, Fan and Li [3].

In this paper, we focus on the case where $\mathscr{P}$ is an *exponential family*. This means that the density functions $p(y, \theta) = \frac{dP_\theta}{dP}(y)$ are of the form $p(y, \theta) = p(y)e^{yC(\theta) - B(\theta)}$. Here, $C(\theta)$ and $B(\theta)$ are some given nondecreasing functions on $\Theta$ and $p(y)$ is some nonnegative function on $\mathbb{R}$.

A natural parametrization for this family is defined by the equality $\mathbf{E}_\theta Y = \int yp(y, \theta)P(dy) = \theta$, for all $\theta \in \Theta$. This condition is useful because the weighted average of observations is a natural unbiased estimate of $\theta$. In what follows, we assume that $\mathscr{P}$ also satisfies the following regularity conditions:

(A1)  $\mathscr{P} = (P_\theta, \theta \in \Theta \subseteq \mathbb{R})$ is an exponential family with a natural parametrization and the functions $B(\cdot)$ and $C(\cdot)$ are continuously differentiable.

(A2)  $\Theta$ is compact and convex and the Fisher information $I(\theta) := \mathbf{E}_\theta |\partial \log p(Y, \theta)/\partial \theta|^2$ satisfies, for some $\varkappa \geq 1$,

$$|I(\theta')/I(\theta'')|^{1/2} \leq \varkappa, \qquad \theta', \theta'' \in \Theta.$$

We illustrate this setup with two examples relevant to the applications we consider below. More examples can be found in Fan, Farmen and Gijbels [5] and Polzehl and Spokoiny [12].

EXAMPLE 2.1 [*Inhomogeneous Bernoulli (binary response) model*].  Let $Z_i = (X_i, Y_i)$ with $X_i \in \mathbb{R}^d$ and $Y_i$ a Bernoulli r.v. with parameter $f(X_i)$, that is, $\mathbf{P}(Y_i = 1 \mid X_i = x) = f(x)$ and $\mathbf{P}(Y_i = 0 \mid X_i = x) = 1 - f(x)$. Such models arise in many econometric applications and are widely used in classification and digital imaging.

EXAMPLE 2.2 (*Inhomogeneous Poisson model*).  Suppose that every $Y_i$ is valued in the set $\mathbb{N}$ of nonnegative integers and $\mathbf{P}(Y_i = k \mid X_i = x) = f^k(x)e^{-f(x)}/k!$, that is, $Y_i$ follows a Poisson distribution with parameter $\theta = f(x)$. This model is commonly used in queueing theory, occurs in positron emission tomography and also serves as an approximation for the density model obtained by a binning procedure.

In the parametric setting with $f(\cdot) \equiv \theta$, the distribution of every "observation" $Y_i$ coincides with $P_\theta$ for some $\theta \in \Theta$ and the parameter $\theta$ can be well estimated using the parametric maximum likelihood method,

$$\widetilde{\theta} = \arg \max_{\theta \in \Theta} \sum_{i=1}^n \log p(Y_i, \theta).$$

In the nonparametric framework, one usually applies the localization idea. In the local constant setting this means that the regression function $f(\cdot)$ can be well approximated by a constant within some neighborhood of every point $x$ in the "feature" space $\mathscr{X}$. This leads to the local model concentrated in some neighborhood of the point $x$.

2.1. *Localization.*  We use *localization by weights* as a general method to describe a local model. Let, for a fixed $x$, a nonnegative weight $w_i = w_i(x) \leq 1$ be assigned to the observation $Y_i$ at $X_i$, $i = 1, \ldots, n$. The weights $w_i(x)$ determine a local model corresponding to the point $x$ in the sense that, when estimating the local parameter $f(x)$, every observation $Y_i$ is taken with weight $w_i(x)$. This leads to the local (weighted) maximum likelihood estimate $\widetilde{\theta} = \widetilde{\theta}(x)$ of $f(x)$,

$$(2.1) \qquad \widetilde{\theta}(x) = \arg\max_{\theta \in \Theta} \sum_{i=1}^{n} w_i(x) \log p(Y_i, \theta).$$

We now mention two possible ways of choosing the weights $w_i(x)$. *Localization by a bandwidth* is defined by weights of the form $w_i(x) = K_{\mathrm{loc}}(\mathbf{l}_i)$ with $\mathbf{l}_i = \rho(x, X_i)/h$, where $h$ is a bandwidth, $\rho(x, X_i)$ is the Euclidean distance between $x$ and the design point $X_i$ and $K_{\mathrm{loc}}$ is a *location kernel. Localization by a window* simply restricts the model to a subset (window) $U = U(x)$ of the design space which depends on $x$, that is, $w_i(x) = \mathbf{1}(X_i \in U(x))$. Observations $Y_i$ with $X_i$ outside the region $U(x)$ are not used for estimating $f(x)$. This kind of localization arises, for example, in the classification with $k$-nearest neighbors method or in the regression tree approach. Sometimes it is convenient to combine these two methods by defining $w_i(x) = K_{\mathrm{loc}}(\mathbf{l}_i)\mathbf{1}(X_i \in U(x))$. One example is given by the boundary-corrected kernels.

We do not assume any special structure for the weights $w_i(x)$, that is, any configuration of weights is allowed. We also denote $W = W(x) = \{w_1(x), \ldots, w_n(x)\}$ and

$$L(W, \theta) = \sum_{i=1}^{n} w_i(x) \log p(Y_i, \theta).$$

To simplify notation, we do not show the dependence of the weights on $x$ explicitly in what follows.

2.2. *Local likelihood estimation for an exponential family model.* If $\mathscr{P} = (P_\theta)$ is an exponential family with the natural parametrization, the local log-likelihood and the local maximum likelihood estimates admit a simple closed form representation. For a given set of weights $W = \mathrm{diag}\{w_1, \ldots, w_n\}$ with $w_i \in [0, 1]$, denote

$$N = \sum_{i=1}^{n} w_i, \qquad S = \sum_{i=1}^{n} w_i Y_i.$$

Note that both sums depend on the location $x$ via the weights $\{w_i\}$.

Lemma 2.1 (Polzehl and Spokoiny [12]).

$$L(W, \theta) = \sum_{i=1}^{n} w_i \log p(Y_i, \theta) = SC(\theta) - NB(\theta) + R,$$

where $R = \sum_{i=1}^{n} w_i \log p(Y_i)$. Moreover,

$$(2.2) \qquad \widetilde{\theta} = S/N = \sum_{i=1}^{n} w_i Y_i \Big/ \sum_{i=1}^{n} w_i$$

and

$$L(W, \widetilde{\theta}, \theta) := L(W, \widetilde{\theta}) - L(W, \theta) = N\mathscr{K}(\widetilde{\theta}, \theta).$$

We now present some exponential inequalities for the "fitted log-likelihood" $L(W, \widetilde{\theta}, \theta)$ which apply in the parametric situation $f(\cdot) \equiv \theta$ for an arbitrary weighting scheme and arbitrary sample size.

Theorem 2.2 (Polzehl and Spokoiny [12]). *Let $W = \{w_i\}$ be a localizing scheme such that $\max_i w_i \leq 1$. If $f(X_i) \equiv \theta^*$ for all $X_i$ with $w_i > 0$, then for any $z > 0$,*

$$\mathbf{P}_{\theta^*}(L(W, \widetilde{\theta}, \theta^*) > z) = \mathbf{P}_{\theta^*}(N\mathscr{K}(\widetilde{\theta}, \theta^*) > z) \leq 2e^{-z}.$$

Remark 2.1. Condition (A2) ensures that the Kullback–Leibler divergence $\mathscr{K}$ satisfies $\mathscr{K}(\theta', \theta^*) \leq I^* |\theta' - \theta^*|^2$ for any point $\theta'$ in a neighborhood of $\theta^*$, where $I^*$ is the maximum of the Fisher information over this neighborhood. Therefore, the result of Theorem 2.2 guarantees that $|\widetilde{\theta} - \theta^*| \leq CN^{-1/2}$ with high probability. Theorem 2.2 can be used to construct the confidence intervals for the parameter $\theta^*$.

Theorem 2.3. *If $\mathfrak{z}_\alpha$ satisfies $2e^{-\mathfrak{z}_\alpha} \leq \alpha$, then*

$$\mathscr{E}_\alpha = \{\theta' : N\mathscr{K}(\widetilde{\theta}, \theta') \leq \mathfrak{z}_\alpha\}$$

*is an $\alpha$-confidence set for the parameter $\theta^*$.*

Theorem 2.2 claims that the estimation loss measured by $\mathscr{K}(\theta', \theta)$ is, with high probability, bounded by $\mathfrak{z}/N$, provided that $\mathfrak{z}$ is sufficiently large. Similarly, one can establish a risk bound for a power loss function.

Theorem 2.4. *Assume* (A1) *and* (A2) *and let $Y_i$ be i.i.d. from $P_{\theta^*}$. Then, for any $r > 0$,*

$$\mathbf{E}_{\theta^*} L^r(\widetilde{\theta}, \theta^*) \equiv N^r \mathbf{E}_{\theta^*} \mathscr{K}^r(\widetilde{\theta}, \theta^*) \leq \tau_r,$$

*where $\tau_r = 2r \int_{\mathfrak{z} \geq 0} \mathfrak{z}^{r-1} e^{-\mathfrak{z}} d\mathfrak{z} = 2r \Gamma(r)$. Moreover, for every $\lambda < 1$,*

$$\mathbf{E}_{\theta^*} \exp\{\lambda L(\widetilde{\theta}, \theta^*)\} \equiv \mathbf{E}_{\theta^*} \exp\{\lambda N\mathscr{K}(\widetilde{\theta}, \theta^*)\} \leq 2(1-\lambda)^{-1}.$$

Proof. By Theorem 2.2,

$$\mathbf{E}_{\theta^*} L^r(\widetilde{\theta}, \theta^*) \leq -\int_{\mathfrak{z} \geq 0} \mathfrak{z}^r \, d\mathbf{P}_{\theta^*}(L(\widetilde{\theta}, \theta^*) > \mathfrak{z})$$

$$\leq r \int_{\mathfrak{z} \geq 0} \mathfrak{z}^{r-1} \mathbf{P}_{\theta^*}(L(\widetilde{\theta}, \theta^*) > \mathfrak{z}) \, d\mathfrak{z} \leq 2r \int_{\mathfrak{z} \geq 0} \mathfrak{z}^{r-1} e^{-\mathfrak{z}} \, d\mathfrak{z}$$

and the first assertion is satisfied. The last assertion is proved similarly. □

2.3. *Risk of estimation in the nonparametric situation.* "*Small modeling bias*" *condition.* This section extends the bound of Theorem 2.2 to the nonparametric situation where the function $f(\cdot)$ is no longer constant, even in the vicinity of the reference point $x$. We, however, suppose that the function $f(\cdot)$ can be well approximated by a constant $\theta$ at all points $X_i$ with positive weights $w_i$. To measure the quality of the approximation, define for every $\theta$

$$(2.3) \qquad \Delta(W, \theta) = \sum_i \delta(\theta, f(X_i)) \mathbf{1}(w_i > 0),$$

where, with $\ell(y, \theta, \theta') = \log \frac{p(y,\theta)}{p(y,\theta')}$,

$$\delta(\theta, \theta') = \log E_\theta e^{-2\ell(Y,\theta,\theta')} = \log E_\theta \frac{p^2(Y,\theta')}{p^2(Y,\theta)}.$$

One can easily check that $\delta(\theta, \theta') \leq I^* |\theta - \theta'|^2$, where $I^* = \max_{\theta'' \in [\theta, \theta']} I(\theta'')$.

THEOREM 2.5. *Let $\mathscr{F}_W$ be a $\sigma$-field generated by the r.v. $Y_i$ for which $w_i > 0$ and let $\Delta(W, \theta) \leq \Delta$. Then, for any random variable $\xi$ measurable with respect to $\mathscr{F}_W$,*

$$\mathbf{E}_{f(\cdot)} \xi \leq (e^\Delta \mathbf{E}_\theta \xi^2)^{1/2}.$$

PROOF. Define $Z_W(\theta) = \exp\{-\sum_i \ell(Y_i, \theta, f(X_i)) \mathbf{1}(w_i > 0)\}$. This value is nothing but the likelihood ratio of the measure $\mathbf{P}_{f(\cdot)}$ with respect to $\mathbf{P}_\theta$ upon restricting to the observations $Y_i$ for which $w_i > 0$. Then, for any $\xi \sim \mathscr{F}_W$, we have $\mathbf{E}_{f(\cdot)} \xi = \mathbf{E}_\theta \xi Z_W(\theta)$. Independence of the $Y_i$'s implies that

$$\log \mathbf{E}_\theta Z_W^2(\theta) = \sum_i \log \mathbf{E}_\theta e^{-2\ell(Y_i, \theta, f(X_i))} \mathbf{1}(w_i > 0)$$

$$= \sum_i \delta(\theta, f(X_i)) \mathbf{1}(w_i > 0) \leq \Delta.$$

The result now follows from the Cauchy–Schwarz inequality $\mathbf{E}_\theta \xi Z_W(\theta) \leq \{\mathbf{E}_\theta \xi^2 \mathbf{E}_\theta Z_W^2(\theta)\}^{1/2}$. □

This result implies that the bound for the risk of estimation $\mathbf{E}_{f(\cdot)}L^r(\widetilde{\theta}, \theta) \equiv N^r \mathbf{E}_{f(\cdot)} \mathscr{K}^r(\widetilde{\theta}, \theta)$ under the parametric hypothesis can be extended to the nonparametric situation, provided that the value $\Delta(W, \theta)$ is sufficiently small.

COROLLARY 2.6. *For any $r > 0$ and any $\lambda < 1$,*

$$N^r \mathbf{E}_{f(\cdot)}|\mathscr{K}(\widetilde{\theta}, \theta)|^r \leq \sqrt{e^{\Delta(W, \theta)} \tau_{2r}},$$

$$N\{\mathbf{E}_{f(\cdot)}|\mathscr{K}(\widetilde{\theta}, \theta)|^r\}^{1/r} \leq \frac{1}{\lambda}\left\{\log \frac{2}{1-\lambda} + \Delta(W, \theta) + 2(r-1)_+\right\}.$$

PROOF. The first bound follows directly from Theorems 2.4 and 2.5. The proof of the second bound uses the fact that for $r > 0$, the function $h(x) = \log^r(x + c_r)$ with $c_r = e^{(r-1)_+}$ is concave on $(0, \infty)$ because

$$h''(x) = \frac{r \log^{r-2}(x + c_r)}{(x + c_r)^2}\{r - 1 - \log(x + c_r)\} \leq 0$$

for $x \geq 0$. With $\zeta = \lambda L(\widetilde{\theta}, \theta)/2$, this implies, by monotonicity of the logarithm function and Jensen's inequality, that $\mathbf{E}_{f(\cdot)}\zeta^r \leq \mathbf{E}_{f(\cdot)}h(e^\zeta) \leq h(\mathbf{E}_{f(\cdot)}e^\zeta)$, hence

$$\mathbf{E}_{f(\cdot)}^{1/r}\zeta^r \leq \log(\mathbf{E}_{f(\cdot)}e^\zeta + c_r) \leq \log \mathbf{E}_{f(\cdot)}e^\zeta + (r-1)_+$$

$$\leq \tfrac{1}{2}\log(e^{\Delta(W, \theta)}\mathbf{E}_\theta e^{2\zeta}) + (r-1)_+$$

and the assertion follows. □

Corollary 2.6 presents two bounds for the risk of estimation in the nonparametric situation which extend the similar parametric bounds by Theorem 2.5. The risk bound in the parametric situation can be interpreted as the bound for the variance of the estimate $\widetilde{\theta}$, while the term $\Delta(W, \theta)$ controls the bias of estimation; see the next section for more details. Both bounds formally apply, regardless of what the "modeling bias" $\Delta(W, \theta)$ is. However, the results are meaningful only if this bias is not too large. The first bound could be preferable for small values of $\Delta(W, \theta)$. However, the multiplicative factor $e^{\Delta(W, \theta)/2}$ makes this bound useless for large $\Delta(W, \theta)$. The advantage of the second bound is that the "modeling bias" enters in additive form.

In the remainder of this section, we briefly comment on relations between the results of Section 2.3 and the usual rate results under smoothness conditions on the function $f(\cdot)$ and the regularity conditions on the design $X_1, \ldots, X_n$. More precisely, we assume that the weights $w_i$ are supported on a ball with radius $h > 0$ and center $x$ and that the function $f(\cdot)$ is smooth within this ball in the sense that for $\theta^* = f(x)$,

$$(2.4) \qquad \delta^{1/2}(\theta^*, f(x + t)) \leq Lh \qquad \forall |t| \leq h.$$

In view of the inequality $\delta(\theta, \theta') \leq I^*|\theta - \theta'|^2$, this condition is equivalent to the usual Lipschitz property. Obviously, (2.4) implies, with $\overline{N} = \sum_i \mathbf{1}(w_i > 0)$, that

$$\Delta(W, \theta^*) \leq L^2 h^2 \overline{N}.$$

Combined with the result of Corollary 2.6, these bounds lead to the following rate results.

THEOREM 2.7. *Assume* (A1) *and* (A2) *and let* $\delta^{1/2}(\theta^*, f(x+t)) \leq Lh$ *for all* $|t| \leq h$. *Select* $h = c(L^2 n)^{-1/(2+d)}$ *for some* $c > 0$ *and let the localizing scheme* $W$ *be such that* $w_i = 0$ *for all* $X_i$ *with* $|X_i - x| > h$, $N := \sum_i w_i \geq \mathfrak{d}_1 n h^d$ *and* $\overline{N} := \sum_i \mathbf{1}(w_i > 0) \leq \mathfrak{d}_2 n h^d$ *with some constants* $\mathfrak{d}_1 < \mathfrak{d}_2$. *Then*

$$\mathbf{E}_{f(\cdot)} |N \mathscr{K}(\widetilde{\theta}, \theta^*)|^{r/2} \leq \{\exp(c^{2+d} \mathfrak{d}_2) \tau_r\}^{1/2}.$$

*Moreover, with* $c_2 = c^{rd/2} \exp(c^{2+d} \mathfrak{d}_2/2) \mathfrak{d}_1^{-r/2}$, *we have*

$$\mathbf{E}_{f(\cdot)} |n^{1/(2+d)} \mathscr{K}(\widetilde{\theta}, \theta^*)|^{r/2} \leq c_2 L^{rd/(2+d)} \tau_r^{1/2}.$$

This corresponds to the classical accuracy of nonparametric estimation for the Lipschitz functions (cf. Fan, Farmen and Gijbels [5]).

**3. Description of the method.** We start by describing the setting under consideration. Let a point of interest $x$ be fixed. The target of estimation is the value $f(x)$ of the regression function at $x$. The local parametric approach described in Section 2 and based on local constant approximation of the regression function in a vicinity of the point $x$ strongly relies on the choice of the local neighborhood or, more generally, of the set of weights $(w_i)$. The problem of selecting such weights and constructing an adaptive (data-driven) estimate is one of the main issues for practical applications and we focus on this problem in this section.

3.1. *Local adaptive estimation. General setup.* For a fixed $x$, we assume that an *ordered* set of localizing schemes $W^{(k)} = (w_i^{(k)})$, for $k = 1, \ldots, K$, is given. The ordering condition means that $w_i^{(k)} \geq w_i^{(k')}$ for all $i$ and all $k > k'$, that is, the degree of locality given by $W_i^{(k)}$ is weakened as $k$ grows. See Section 3.3 for some examples. For the popular example of kernel weights $w_i^{(k)} = K((X_i - x)/h_k)$, this condition means that the bandwidth $h_k$ grows with $k$. Also, let $\{\widetilde{\theta}^{(k)}, k = 1, \ldots, K\}$ be the corresponding set of local likelihood estimates for $\theta = f(x)$,

$$\widetilde{\theta}^{(k)}(x) = \arg\max_{\theta \in \Theta} L(W, \theta) = \sum_{i=1}^{n} w_i^{(k)} Y_i \Big/ \sum_{i=1}^{n} w_i^{(k)}.$$

Due to Theorem 2.2, the value $1/N_k$ can be used to measure the variability of the estimate $\widetilde{\theta}^{(k)}$. The ordering condition means, in particular, that $N_k$ grows and hence the variability of $\widetilde{\theta}^{(k)}$ decreases with $k$.

Given the estimates $\widetilde{\theta}^{(k)}$, we consider the larger class of their convex combinations,

$$\widehat{\theta} = \alpha_1 \widetilde{\theta}^{(1)} + \cdots + \alpha_K \widetilde{\theta}^{(K)}, \qquad \alpha_1 + \cdots + \alpha_K = 1, \qquad \alpha_k \geq 0,$$

where the mixing coefficients $\alpha_k$ may depend on the point $x$. We aim to construct a new estimate $\widehat{\theta}$ in this class which performs at least as well as the best one in the original family $\{\widetilde{\theta}^{(k)}\}$.

3.2. *Stagewise aggregation procedure.* The adaptive estimate $\widehat{\theta}$ of $\theta = f(x)$ is computed sequentially via the following algorithm:

1. *Initialization*: $\widehat{\theta}^{(1)} = \widetilde{\theta}^{(1)}$.
2. *Stagewise aggregation*: For $k = 2, \ldots, K$,

$$\widehat{\theta}^{(k)} := \gamma_k \widetilde{\theta}^{(k)} + (1 - \gamma_k) \widehat{\theta}^{(k-1)},$$

with the mixing parameter $\gamma_k$ being defined for some $\mathfrak{z}_k > 0$ and a kernel $K_{\mathrm{ag}}(\cdot)$ as

$$\gamma_k = K_{\mathrm{ag}}(\mathbf{m}^{(k)}/\mathfrak{z}_k), \qquad \mathbf{m}^{(k)} := N_k \mathscr{K}(\widetilde{\theta}^{(k)}, \widehat{\theta}^{(k-1)}).$$

3. *Loop*: If $k < K$, then increase $k$ by one and continue with Step 2. Otherwise, terminate and set $\widehat{\theta} = \widehat{\theta}^{(K)}$.

The idea behind the procedure is quite simple. We start with the first estimate $\widetilde{\theta}^{(1)}$ which has the smallest degree of locality but the largest variability, of order $1/N_1$. Next, we consider estimates with larger values $N_k$. Every current estimate $\widetilde{\theta}^{(k)}$ is compared with the previously constructed estimate $\widehat{\theta}^{(k-1)}$. If the difference is not significant, then the new estimate $\widehat{\theta}^{(k)}$ basically coincides with $\widetilde{\theta}^{(k)}$. Otherwise, the procedure essentially keeps the previous value $\widehat{\theta}^{(k-1)}$. For measuring the difference between the estimates $\widetilde{\theta}^{(k)}$ and $\widehat{\theta}^{(k-1)}$, we use $\mathbf{m}^{(k)} := N_k \mathscr{K}(\widetilde{\theta}^{(k)}, \widehat{\theta}^{(k-1)})$, which is motivated by the results of Theorems 2.2 and 2.3. In particular, a large value of $\mathbf{m}^{(k)}$ means that $\widehat{\theta}^{(k-1)}$ does not belong to the confidence set corresponding to $\widetilde{\theta}^{(k)}$ and hence indicates a significant difference between these two estimates. To quantify this significance, the procedure utilizes the parameters (critical values) $\mathfrak{z}_k$. Their choice is discussed in Section 3.3.1.

REMARK 3.1. If $K_{\mathrm{ag}}(\cdot)$ is the uniform kernel on $[0, 1]$, then $\gamma_k$ is either zero or one, depending on the value of $\mathbf{m}^{(k)}$. This yields, by induction arguments, that the final estimate coincides with one of the "weak" estimates $\widetilde{\theta}^{(k)}$. In this case, our method can be considered a pointwise model selection method.

If the kernel $K_{\mathrm{ag}}$ is such that $K_{\mathrm{ag}}(t) = 1$ for $t \le b$ with some positive $b$, then the small values of the "test statistic" $\mathbf{m}^{(k)}$ lead to the aggregated estimate $\widehat{\theta}^{(k)} = \widetilde{\theta}^{(k)}$. This is an important feature of the procedure which will be used in our implementation and theoretical study.

3.3. *Parameter choice and implementation details.* The implementation of the SSA procedure requires the fixing of a sequence of local likelihood estimates, the kernel $K_{\mathrm{ag}}$ and the parameters $\mathfrak{z}_k$. The next section gives some examples of how the set of localizing schemes $W^{(k)}$ can be selected. The only important parameters of the method are "critical values" $\mathfrak{z}_k$ which normalize the "test statistics" $\mathbf{m}^{(k)}$. Section 3.3.1 describes in detail how they can be selected in practice.

The kernel $K_{\mathrm{ag}}$ should satisfy $0 \le K_{\mathrm{ag}}(t) \le 1$, be monotone decreasing and have support on $[0,1]$. Further, there should be a positive number $b$ such that $K_{\mathrm{ag}}(t) = 1$ for $t \le b$. Our default choice is a piecewise linear kernel with $b = 1/6$ and $K_{\mathrm{ag}}(t) = (1 - (t - b)_+)_+$. Our numerical results (not shown here) indicate that the particular choice of kernel $K_{\mathrm{ag}}$ has only a minor effect on the final results.

3.3.1. *Choice of the parameters $\mathfrak{z}_k$.* The "critical values" $\mathfrak{z}_k$ define the level of significance for the test statistics $\mathbf{m}^{(k)}$. A proper choice of these parameters is crucial for the performance of the procedure. In this section, we propose one general approach for selecting them which is similar to the bootstrap idea in the hypothesis testing problem. Namely, we select these values to provide the prescribed performance of the procedure in the parametric situation (under the null hypothesis). For every step $k$, we require that the estimate $\widehat{\theta}^{(k)}$ is sufficiently close to the "oracle" estimate $\widetilde{\theta}^{(k)}$ in the parametric situation $f(\cdot) \equiv \theta$, in the sense that

$$(3.1) \qquad \sup_{\theta^* \in \Theta} \mathbf{E}_{\theta^*} |N_k \mathscr{K}(\widetilde{\theta}^{(k)}, \widehat{\theta}^{(k)})|^r \le \alpha \tau_r$$

for all $k = 2, \ldots, K$ with $\tau_r$ from Theorem 2.4. In some cases, the risk $\mathbf{E}_{\theta^*} |N_k \mathscr{K}(\widetilde{\theta}^{(k)}, \widehat{\theta}^{(k)})|^r$ does not depend on $\theta^*$. This is the case, for example, when $\theta$ is a shift or scale parameter, as for Gaussian shift, exponential and volatility families. It then suffices to check (3.1) for any single point $\theta^*$. In the general situation, the risk $\mathbf{E}_{\theta^*} |N_k \mathscr{K}(\widetilde{\theta}^{(k)}, \widehat{\theta}^{(k)})|^r$ depends on the parameter value $\theta^*$. However, our numerical results (not reported here) indicate that this dependence is minor and usually it suffices to check these conditions for one parameter $\theta^*$. In particular, for the Bernoulli model considered in Section 4 we recommend only checking condition (3.1) for the

"least favorable" value $\theta^* = 1/2$ corresponding to the largest variance of the estimate $\widetilde{\theta}$.

The values $\alpha$ and $r$ in (3.1) are two global parameters. The role of $\alpha$ is similar to the level of the test in the hypothesis testing problem, while $r$ describes the power of the loss function. A specific choice is subjective and depends on the particular application in question. Taking a large $r$ and small $\alpha$ would result in an increase of the critical values and therefore improve the performance of the method in the parametric situation, with the cost of some loss of sensitivity to parameter changes. Theorem 5.1 presents some upper bounds for the critical values $\mathfrak{z}_k$ as functions of $\alpha$ and $r$ in the form $a_0 + a_1 \log \alpha^{-1} + a_2 r(K - k)$ with some coefficients $a_0$, $a_1$ and $a_2$. We see that these bounds depend linearly on $r$ and $\log \alpha^{-1}$. For our applications to classification, we apply a relatively small value, $r = 1/2$, because the misclassification error corresponds to the bounded loss function. We also apply $\alpha = 1$, although other values in the range $[0.5, 1]$ lead to very similar results. Note that in general, the parameters $\mathfrak{z}_k$ thus defined depend on the model considered, design $X_1, \ldots, X_n$ and localizing schemes $W^{(1)}, \ldots, W^{(K)}$, which, in turn, can differ from point to point. Therefore, an implementation of the suggested rule would require separate computation of the parameters for every point of estimation. However, in many situations, for example, for the regular design, this variation from point to point is negligible and a universal set of parameters can be used. It is only important that conditions (3.1) are satisfied for all the points.

3.3.2. *Simplified parameter choice.* Proposal (3.1) is not constructive: we have only $K - 1$ conditions for choosing $K - 1$ parameters. Here, we present a simplified procedure which is rather easy to implement and is based on Monte Carlo simulations. It suggests first identifying the last value $\mathfrak{z}_K$ using the reduced aggregation procedure with only two estimates $\widetilde{\theta}^{(K-1)}$ and $\widetilde{\theta}^{(K)}$:

$$\sup_{\theta^* \in \Theta} \mathbf{E}_{\theta^*} |N_K \mathscr{K}(\widetilde{\theta}^{(K)}, \widehat{\theta}(\mathfrak{z}_K))|^r \leq \alpha \tau_r / (K - 1),$$

where $\widehat{\theta}(\mathfrak{z}_K) = \gamma \widetilde{\theta}^{(K)} + (1 - \gamma) \widetilde{\theta}^{(K-1)}$, $\gamma = K_{\mathrm{ag}}(\mathbf{m}/\mathfrak{z}_K)$ and $\mathbf{m} = N_K \mathscr{K}(\widetilde{\theta}^{(K)}, \widetilde{\theta}^{(K-1)})$. The other values $\mathfrak{z}_k$ are found in the form $\mathfrak{z}_k = \mathfrak{z}_K + \iota(K - k)$ to ensure (3.1). This suggestion is justified by Theorem 5.1 from Section 5.1.

3.3.3. *Examples of sequences of local likelihood estimates.* This section presents some examples and recommendations for the choice of the localizing schemes $W^{(k)}$ which we also use in our simulation study. Note, however, that the choice of $W^{(k)}$'s is not part of the SSA procedure. The procedure applies with any choice under some rather mild growth conditions.

Below, we assume that the design $X_1, \ldots, X_n$ is supported on the unit cube $[-1, 1]^d$. This condition can be easily met by rescaling the design components. We mention two approaches to choosing the localizing scheme which are usually used in applications. One is based on a given sequence of bandwidths, and the other is based on the nearest neighbor structure of the design. In both situations, we assume that a *location* kernel $K_{\mathrm{loc}}$ is a nonnegative function on the unit cube $[-1, 1]^d$. In general, we only assume that this kernel is decreasing along any radial line, that is, $K_{\mathrm{loc}}(\rho x) \geq K_{\mathrm{loc}}(x)$ for any $x \in [-1, 1]^d$ and $\rho \leq 1$, and $K_{\mathrm{loc}}(x) = 0$ for $|x| \geq 1$. In most applications, one applies an isotropic kernel $K_{\mathrm{loc}}$ which depends only on the norm of $x$. The recommended choice is the Epanechnikov kernel $K_{\mathrm{loc}}(x) = (1 - |x|^2)_+$.

*Bandwidth-based localizing schemes.* This approach is recommended for the univariate or bivariate equidistant design. Let $\{h_k\}_{k=1}^K$ be a finite set of *bandwidth candidates.* We assume that this set is ordered, that is, $h_1 < h_2 < \cdots < h_K$. Every such bandwidth determines the collection of kernel weights $w_i^{(k)} = K_{\mathrm{loc}}((X_i - x)/h_k)$, $i = 1, \ldots, n$. This definition assumes that the same localizing bandwidth is applied for all directions in the feature space. In all of the examples below, we apply a geometrically increasing sequence of "bandwidths" $h_k$, that is, $h_{k+1} = ah_k$ for some $a > 1$. This sequence is uniquely determined by the starting value $h_1$, the factor $a$ and the total number $K$ of local schemes. The recommended choice of $a$ is $(1.25)^{1/d}$, although our numerical results (not reported here) indicate no significant change in the results when any other value of $a$ in the range 1.1 to 1.3 is used. The value $h_1$ is to be selected in such a way that the starting estimate $\widetilde{\theta}^{(1)}$ is well defined for all points of estimation. In the case of a local constant approximation, this value can be taken very small because even one point can be sufficient for preliminary estimation. In the case of a regular design, the value $h_1$ is of order $n^{-1/d}$. The number $K$ of local schemes $W^{(k)}$ or, equivalently, of "weak" estimates $\widetilde{\theta}^{(k)}$ is largely determined by the values $h_1$ and $a$, in such a way that $h_K = h_1 a^{K-1}$ is approximately one, that is, the last estimate behaves like a global parametric estimate from the whole sample. The formula $K = a \log(h_K/h_1)$ suggests that $K$ is at most logarithmic in the sample size $n$.

*k-NN-based local schemes.* If the design is irregular or the design space is high-dimensional ($d > 2$), then it is useful to apply the local schemes based on the $k$-nearest neighbor structure of the design. For this approach, an increasing sequence $\{N_k\}$ of integers must be fixed. For a fixed $x$ and every $k \geq 1$, the bandwidth $h_k$ is the minimal one for which the ball of radius $h_k$ contains at least $N_k$ design points. The weights are again defined by $w_i^{(k)} = K_{\mathrm{loc}}((X_i - x)/h_k)$. The sequence $\{N_k\}$ is selected similarly to the sequence $\{h_k\}$ in the bandwidth-based approach. One starts with a fixed $N_1$ and then multiplies it at every step by some factor $a > 1$: $N_{k+1} = aN_k$. The number of steps $K$ is such that $N_K$ is of order $n$.

One can easily check that the kernel- and $k$-NN-based local schemes coincide in the case of univariate regular design.

**4. Application to classification.** One observes a training sample $(X_i, Y_i)$, $i = 1, \ldots, n$, with $X_i$ valued in a Euclidean space $\mathscr{X} = \mathbb{R}^d$ and with known class assignment $Y_i \in \{0, 1\}$. Our objective is to construct a discrimination rule assigning every point $x \in \mathscr{X}$ to one of the two classes. The classification problem can be naturally treated in the context of a binary response model. It is assumed that each observation $Y_i$ at $X_i$ is a Bernoulli r.v. with the parameter $\theta_i = f(X_i)$, that is, $\mathbf{P}(Y_i = 0 | X_i) = 1 - f(X_i)$ and $\mathbf{P}(Y_i = 1 | X_i) = f(X_i)$. The "ideal" Bayes discrimination rule is $\mathbf{1}(f(x) \geq 1/2)$. Since the function $f(x)$ is usually unknown, it is replaced by its estimate $\widehat{\theta}$. If the distribution of $X_i$ within the class $k$ has density $p_k$, then

$$\theta_i = \pi_1 p_1(X_i) / (\pi_0 p_0(X_i) + \pi_1 p_1(X_i)),$$

where $\pi_k$ is the prior probability of $k$th population, $k = 0, 1$.

Nonparametric methods of estimating the function $\theta$ are typically based on local averaging. Two typical examples are given by the $k$-nearest neighbor ($k$-NN) estimate and the kernel estimate. For a given $k$ and every point $x$ in $\mathscr{X}$, denote by $\mathscr{D}_k(x)$ the subset of the design $X_1, \ldots, X_n$ containing the $k$ nearest neighbors of $x$. Then the $k$-NN estimate of $f(x)$ is defined by averaging the observations $Y_i$ over $\mathscr{D}_k(x)$,

$$\widetilde{\theta}^{(k)}(x) = k^{-1} \sum_{X_i \in \mathscr{D}_k(x)} Y_i.$$

The definition of the kernel estimate of $f(x)$ involves a univariate kernel function $K(\cdot)$ and the bandwidth $h$,

$$\widetilde{\theta}^{(h)}(x) = \sum_{i=1}^{n} K\left(\frac{\rho(x, X_i)}{h}\right) Y_i \Big/ \sum_{i=1}^{n} K\left(\frac{\rho(x, X_i)}{h}\right).$$

Both methods require the choice of a smoothing parameter (the value $k$ for $k$-NN and the bandwidth $h$ for the kernel estimate).

EXAMPLE 4.1. In this example, we consider the binary classification problem with the corresponding class densities $p_0(x)$ and $p_1(x)$ given by two-component normal mixtures

$$p_0(x) = 0.2\phi(x; (-1, 0), 0.5\mathbf{I}_2) + 0.8\phi(x; (1, 0), 0.5\mathbf{I}_2),$$

$$p_1(x) = 0.5\phi(x; (0, 1), 0.5\mathbf{I}_2) + 0.5\phi(x; (0, -1), 0.5\mathbf{I}_2),$$

where $\phi(\cdot; \mu, \Sigma)$ is the density of the multivariate normal distribution with mean vector $\mu$ and covariance matrix $\Sigma$ and $\mathbf{I}_2$ is the $2 \times 2$ unit matrix.

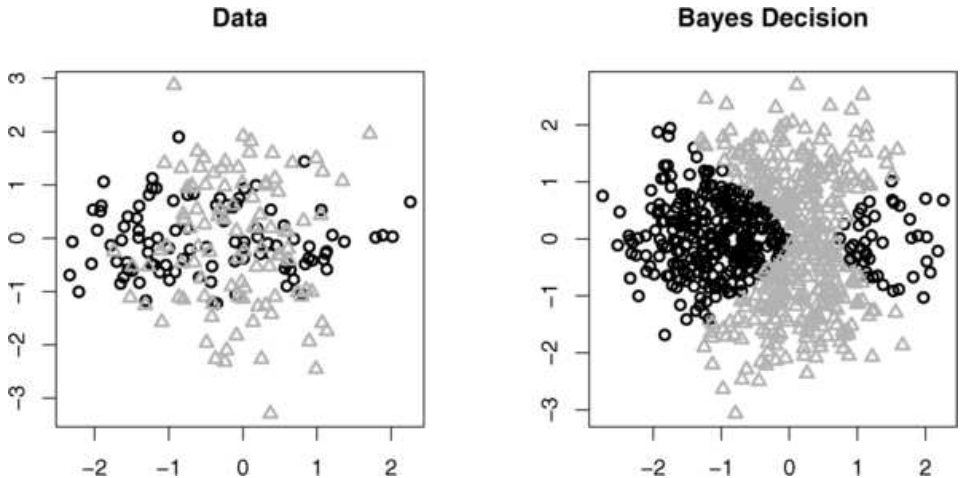

Fig. 1. *A sample from the binary response model with the normal mixture class densities (left) and results of applying the Bayes discrimination rule to this model (right).*

Figure 1 shows one typical realization of the training sample with 100 observations in each class (left) and the optimal Bayes classification for a testing sample with 1000 observations in each class (right). First, in order to illustrate the "oracle" property of the SSA, we compute the pointwise misclassification errors for all weak estimates and SSA estimates at four boundary points. Figure 2 is obtained using a training sample of size 400, $k$-NN weighting scheme with $N_1 = 5$, $N_K = 300$, $K = 30$ and $\alpha = 0.5$. Further, we have carried out 500 simulation runs, each time generating 100 training points and 100 testing points. The rates of misclassification on testing sets have been averaged thereafter to give the *mean misclassification error*, shown as a dotted reference line in Figure 3. We note here that the critical values

$$\mathfrak{z}_k = 0.0031 + 0.007 * (K - k), \qquad k = 1, \ldots, K,$$

have been computed only once for one design realization and least favorable parameter value $\theta^* = 0.5$, then used in all runs. The same strategy is used in other examples. Next, two "weak" classification methods, $k$-NN and kernel classifiers, with varying smoothing parameters, are applied to the same data set. Figure 3 (top) shows the dependence of the misclassification error on the bandwidth for kernel classifiers and on the number of nearest neighbors for the $k$-NN classifier.

One can observe that a careful choice of the smoothing parameter is crucial for getting a reasonable quality of classification. A wrong choice leads to a significant increase of the misclassification rate, especially for the kernel classifiers. At the same time, the optimal choice can lead to a reasonable

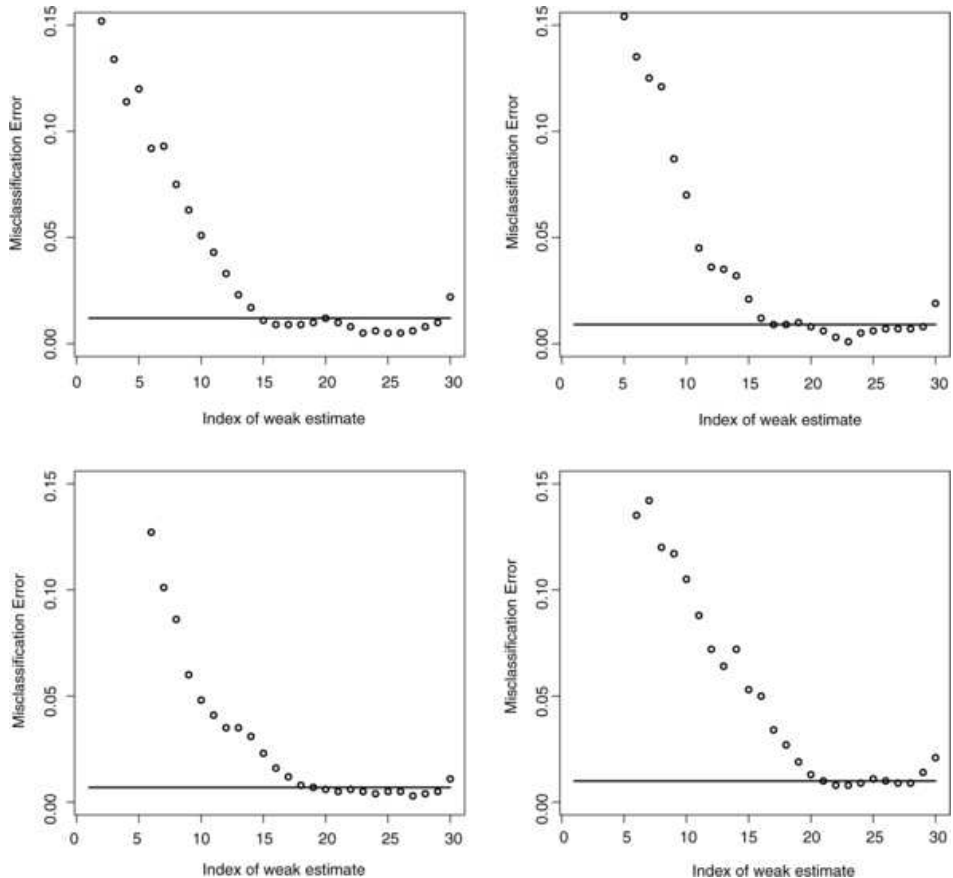

Fig. 2. *Pointwise misclassification errors (black dots) at four points for all weak estimates used in Example 4.1. Solid reference lines correspond to SSA misclassification errors.*

quality of classification, only slightly worse than that of the Bayes decision rule.

EXAMPLE 4.2.   We now consider Example 4.1 with eight additional independent $\mathcal{N}(0,1)$-distributed nuisance components. So, now $X_i = (X_i^1, \dots, X_i^{10})$, where

$$(X_i^1, X_i^2) \sim p_{\mathrm{class}(i)}, \qquad (X_i^3, \dots, X_i^{10}) \sim \mathcal{N}((\underbrace{0, \dots, 0}_{8}), \mathbf{I}_8).$$

The SSA procedure is now implemented, again using $k$-NN weights with the number of nearest neighbors exponentially increasing from 5 to 100. The results are shown in the bottom row of Figure 3. We again observe that the quality of both standard classifiers depends significantly on the choice of the smoothing parameters. In the high-dimensional situation considered, even under the optimal choice, the quality of the dimension-independent Bayes classifier is not attained. However, the SSA procedure again performs nearly as well as the best k-NN or kernel classifier.

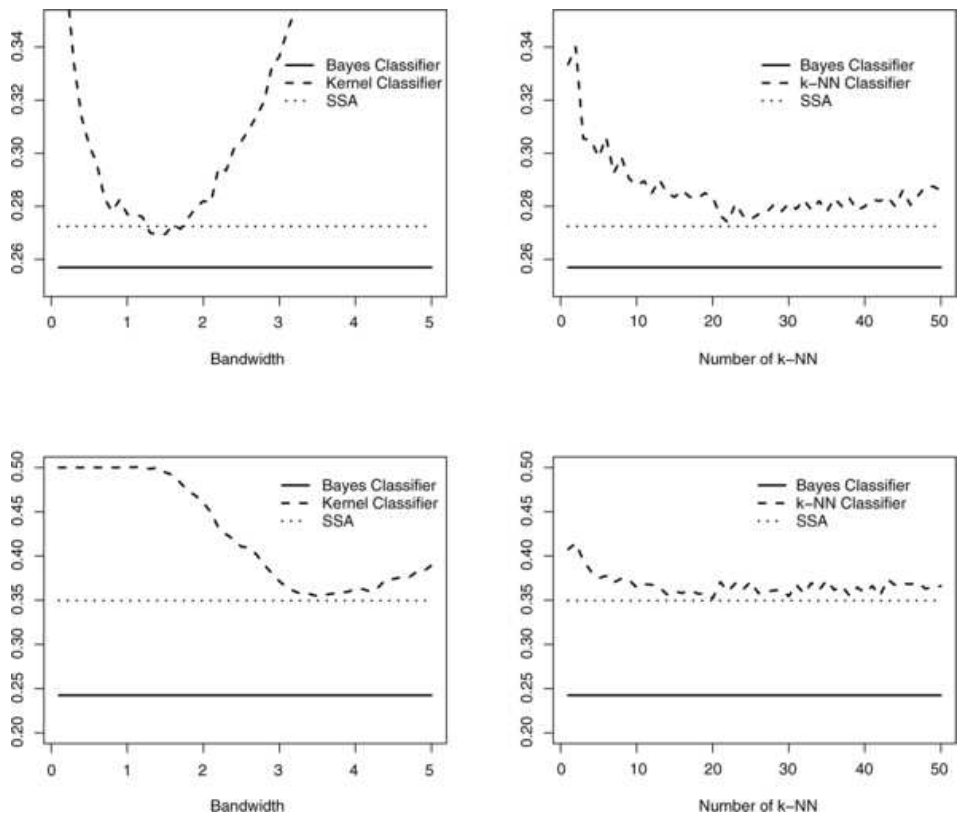

Fig. 3. *Misclassification errors as functions of the main smoothing parameter for k-NN (right) and kernel (left) classifiers. SSA and Bayes misclassification errors are given as reference lines.* Top: *Example 4.1 (dimension 2).* Bottom: *Example 4.2 (dimension 10).*

EXAMPLE 4.3 (*BUPA liver disorders*). We consider the dataset sampled by BUPA Medical Research Ltd. It consists of seven variables and 345 observed vectors. The subjects are single male individuals. The first five variables are measurements taken by blood tests that are thought to be sensitive to liver disorders and which might arise from excessive alcohol consumption. The sixth variable is a sort of selector variable. The seventh variable is the label indicating the class identity. Among all the observations, there are 145 people who belong to the liver-disorder group (corresponding to selector number 2) and 200 people who belong to the liver-normal group. The BUPA liver disorder data set is notoriously difficult to classify, with usual error rates at about 30%. We apply SSA, $k$-NN and kernel classifiers to tackle this problem. In the SSA procedure, the $k$-NN weighting scheme was employed with the number of $k$-NN ranging from 2 to 100. Figure 4 shows the corresponding leave-one-out cross-validation errors for the above methods. One can see that the SSA method is uniformly better than kernel or $k$-NN classifiers.

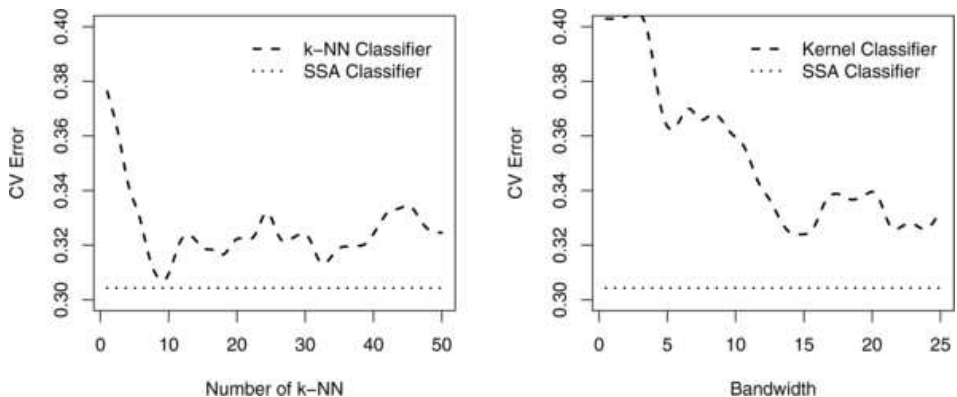

Fig. 4. *One-leave-out cross-validation errors as functions of the main smoothing parameters for k-NN (left) and kernel (right) classifiers. The dotted line describes the error of the SSA classifier.*

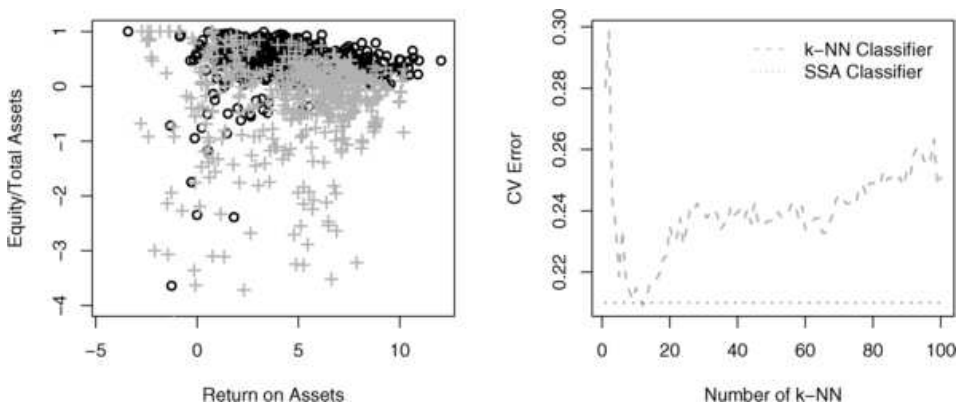

Fig. 5. Left: *default events (crosses indicate defaulted firms and circles indicate operating ones) are shown depending on the two characteristics of a firm.* Right: *leave-one-out cross-validation error for the k-NN classifier as a function of the number of nearest neighbors. The CV error for the SSA classifier is given as a reference line.*

Example 4.4 (*Bankruptcy data*).   The data set from the Compustat repository contains statistics concerning bankruptcies (defaults) in the pri-vate sector of the U.S. economy during the period 2000–2005. There are
14 explanatory variables including different financial ratios, industry indica-
tors and so on. First, a preliminary analysis is conducted and the two most
informative variables [equity/total assets ratio and net income/total assets
ratio (profitability)] are selected. The projection of the default statistics onto
the corresponding plane is shown in Figure 5. Further, the performance of
the SSA procedure is compared to the performance of the $k$-NN classifier
with different numbers of nearest neighbors. Namely, leave-one-out cross-
validation errors are computed for both SSA and $k$-NN classification meth-
ods and the latter one is presented in Figure 5 as a function of the number
of nearest neighbors. Again, as in previous examples, the quality of classi-
fication strongly depends on the choice of the parameter $k$. The adaptive
SSA procedure provides the performance corresponding to the best possible
choice of this parameter.

**5. Some theoretical properties of the SSA method.** This section dis-
cusses some important theoretical properties of the proposed aggregating
procedure. In particular, we establish the "oracle" result which claims that
the aggregated estimate is, up to a log factor, as good as the best one among
the considered family $\{\widetilde{\theta}^{(k)}\}$ of local constant estimates.

The majority of the results in the modern statistical literature are stated
as asymptotic rate results. It is, however, well known that the rate optimal-
ity of an estimation procedure does not automatically imply its good finite
sample properties and cannot be used for comparing different procedures.
Also, the rate results are almost useless for selecting the parameters of the
procedure. In our theoretical study, we apply another approach which aims
to link parametric and nonparametric inference with the focus on the adap-
tive behavior of the proposed method. This means, in particular, that the
SSA procedure attains parametric accuracy if the parametric assumption is
satisfied. In the general situation, the procedure attains (up to an unavoid-
able price for adaptation) the quality corresponding to the best possible
local parametric approximation for the underlying model near the point of
interest.

The "oracle" result is, in turn, a consequence of two important properties
of the aggregated estimate $\widehat{\theta}$: "propagation" and "stability." "Propagation"
can be viewed as the oracle result in the parametric situation with $f(\cdot) \equiv \theta^*$.
In this case, the oracle choice would be the estimate with the largest value
$N_k$, that is, the last estimate $\widetilde{\theta}^{(K)}$ in the family $\{\widetilde{\theta}^{(k)}\}$. The "propagation"
property means that at every step $k$ of the procedure, the "aggregated" es-
timate $\widehat{\theta}^{(k)}$ is close to the "oracle" estimate $\widetilde{\theta}^{(k)}$. In other words, the "propa-
gation" property ensures that at every step, the degree of locality is relaxed
and the local model applied for estimation is extended to a larger neigh-
borhood described by the weights $W^{(k)}$. The "propagation" property can be naturally extended to a nearly parametric case when $\Delta(W^{(k)}, \theta)$ is small for some fixed $\theta$ and all $k \leq k^*$. The "propagation" feature of the procedure ensures that the quality of estimation improves and confidence bounds for $\widehat{\theta}^{(k)}$ become tighter as the number of iterations increases, provided that the "small modeling bias" condition still holds. Finally, the "stability" property ensures that the quality gained in the "propagation" stage will be maintained for the final estimate.

Our theoretical study is carried out under assumptions (A1) and (A2) on the parametric family $\mathscr{P}$. Additionally, we impose the following assumption on the sequence of localizing schemes $W^{(k)}$ which was already mentioned in Section 3.

(A3) The set $W^{(k)}$ is ordered, in the sense that $w_i^{(k)} \geq w_i^{(k')}$ for all $i$ and all $k > k'$. Moreover, for some constants $u_0, u$ with $0 < u_0 \leq u < 1$, values $N_k = \sum_{j=1}^{n} w_i^{(k)}$ satisfy, for every $2 \leq k \leq K$,

$$u_0 \leq N_{k-1}/N_k \leq u.$$

5.1. *Behavior in the parametric situation.* First, we consider the homogeneous situation with the constant parameter value $f(x) = \theta^*$. Our first result claims that in this situation, under assumption (A3), the parameters $\mathfrak{z}_k$ can be chosen in the form $\mathfrak{z}_k = \mathfrak{z}_K + \iota(K - k)$ in order to satisfy the "propagation" condition (3.1). The proof is given in the Appendix.

THEOREM 5.1. *Assume* (A1), (A2) *and* (A3). *Let* $f(X_i) = \theta^*$ *for all $i$. Then there are three constants $a_0$, $a_1$ and $a_2$, depending on $r$, $u_0$ and $u$ only, such that the choice*

$$\mathfrak{z}_k = a_0 + a_1 \log \alpha^{-1} + a_2 r \log N_k$$

*ensures* (3.1) *for all $k \leq K$. In particular,* $\mathbf{E}_{\theta^*} |N_K \mathscr{K}(\widetilde{\theta}^{(K)}, \widehat{\theta})|^r \leq \alpha \tau_r$.

5.2. *"Propagation" under "small modeling bias."* We now extend the "propagation" result to the situation where the parametric assumption is no longer fulfilled, but the deviation from the parametric structure within the local model under consideration is sufficiently small. This deviation can be measured for the localizing scheme $W^{(k)}$ by $\Delta(W^{(k)}, \theta)$ from (2.3).

We suppose that there is a number $k^*$ such that the modeling bias $\Delta(W^{(k)}, \theta)$ is small for some $\theta$ and all $k \leq k^*$. Consider the corresponding estimate $\widehat{\theta}^{(k^*)}$ obtained after the first $k^*$ steps of the algorithm. Theorem 2.5 implies, in this situation, the following result.

THEOREM 5.2. *Assume* (A1), (A2) *and* (A3). *Let $\theta$ and $k^*$ be such that $\Delta(W^{(k)}, \theta) \leq \Delta$ for some $\Delta \geq 0$ and all $k \leq k^*$. Then*

$$\mathbf{E}_{f(\cdot)} |N_{k^*} \mathscr{K}(\widetilde{\theta}^{(k^*)}, \widehat{\theta}^{(k^*)})|^{r/2} \leq \sqrt{\alpha \tau_r e^{\Delta}},$$

$$\mathbf{E}_{f(\cdot)}|N_{k^*}\mathscr{K}(\widetilde{\theta}^{(k^*)}, \theta)|^{r/2} \leq \sqrt{\tau_r e^\Delta}.$$

5.3. *"Stability after propagation" and "oracle" results.* Due to the "propagation" result, the procedure performs well, provided the "small modeling bias" condition $\Delta(W^{(k)}, \theta) \leq \Delta$ is satisfied. To establish the accuracy result for the final estimate $\widehat{\theta}$, we have to check that the aggregated estimate $\widehat{\theta}^{(k)}$ does not vary much at the steps "after propagation" when the divergence $\Delta(W^{(k)}, \theta)$ from the parametric model becomes large.

THEOREM 5.3. *Under* (A1), (A2) *and* (A3), *for every* $k \leq K$, *we have*

$$(5.1) \qquad\qquad N_k \mathscr{K}(\widehat{\theta}^{(k)}, \widehat{\theta}^{(k-1)}) \leq \mathfrak{z}_k.$$

*Moreover, under* (A3), *for every* $k'$ *with* $k < k' \leq K$, *we have*

$$(5.2) \qquad\qquad N_k \mathscr{K}(\widehat{\theta}^{(k')}, \widehat{\theta}^{(k)}) \leq \varkappa^2 c_u^2 \, \overline{\mathfrak{z}}_k$$

*with* $c_u = (u^{-1/2} - 1)^{-1}$ *and* $\overline{\mathfrak{z}}_k = \max_{l \geq k} \mathfrak{z}_l$.

REMARK 5.1. An interesting feature of this result is that it is satisfied with probability one, that is, the control of stability not only "works" with high probability, it always applies. This property follows directly from the construction of the procedure.

PROOF OF THEOREM 5.3. (The convexity of the Kullback–Leibler divergence $\mathscr{K}(u, v)$) with respect to the first argument implies

$$\mathscr{K}(\widehat{\theta}^{(k)}, \widehat{\theta}^{(k-1)}) \leq \gamma_k \mathscr{K}(\widetilde{\theta}^{(k)}, \widehat{\theta}^{(k-1)}).$$

If $\mathscr{K}(\widetilde{\theta}^{(k)}, \widehat{\theta}^{(k-1)}) \geq \mathfrak{z}_k/N_k$, then $\gamma_k = 0$ and (5.1) follows. Now, assumption (A2) and Lemma A.1 yield

$$\mathscr{K}^{1/2}(\widehat{\theta}^{(k')}, \widehat{\theta}^{(k)}) \leq \varkappa \sum_{l=k+1}^{k'} \mathscr{K}^{1/2}(\widehat{\theta}^{(l)}, \widehat{\theta}^{(l-1)}) \leq \varkappa \sum_{l=k+1}^{k'} (\mathfrak{z}_l/N_l)^{1/2}.$$

The use of assumption (A3) leads to the bound

$$\mathscr{K}^{1/2}(\widehat{\theta}^{(k')}, \widehat{\theta}^{(k)}) \leq \varkappa(\overline{\mathfrak{z}}_k/N_k)^{1/2} \sum_{l=k+1}^{k'} u^{(l-k)/2} \leq \varkappa\sqrt{u}(1 - \sqrt{u})^{-1}(\overline{\mathfrak{z}}_k/N_k)^{1/2},$$

which proves (5.2).  □

A combination of the "propagation" and "stability" statements implies the main result concerning the properties of the adaptive estimate $\widehat{\theta}$.

THEOREM 5.4. *Assume* (A1), (A2) *and* (A3). *Let $k^*$ be a "good" choice, in the sense that*

$$\max_{k \le k^*} \Delta(W^{(k)}, \theta) \le \Delta$$

*for some $\theta$ and some value $\Delta$. Then*

$$\mathbf{E}_{f(\cdot)} |N_{k^*} \mathscr{K}(\widetilde{\theta}^{(k^*)}, \widehat{\theta})|^{r/2} \le 2^{(r-1)_+} \varkappa^r \{ \sqrt{\alpha \tau_r e^\Delta} + (c_u^2 \overline{\mathfrak{z}}_{k^*})^{r/2} \},$$

*where $c_u$ is the constant from Theorem* 5.3.

We also present a corollary of the "oracle" result concerning the risk of the adaptive estimate $\widehat{\theta}$ for the special case where $r = 1$. Other values of $r$ can be considered as well: one only has to update the constants depending on $r$. We also assume that $\alpha \le 1$.

COROLLARY 5.5. *Let $\max_{k \le k^*} \Delta(W^{(k)}, \theta) \le \Delta$ for some $\theta$ and some $\Delta$. Then*

$$N_{k^*}^{1/2} \mathbf{E}_{f(\cdot)} \mathscr{K}^{1/2}(\widehat{\theta}, \theta) \le \varkappa (2\sqrt{\tau_1 e^\Delta} + \sqrt{c_u^2 \overline{\mathfrak{z}}_{k^*}}).$$

PROOF. Simply observe that by Lemma A.1,

$$\mathscr{K}^{1/2}(\widehat{\theta}, \theta)$$

$$\le \varkappa \left\{ \mathscr{K}^{1/2}(\widetilde{\theta}^{(k^*)}, \theta) + \mathscr{K}^{1/2}(\widetilde{\theta}^{(k^*)}, \widehat{\theta}^{(k^*)}) + \sum_{l=k^*+1}^{\widehat{k}} \mathscr{K}^{1/2}(\widehat{\theta}^{(l)}, \widehat{\theta}^{(l-1)}) \right\}$$

and follow the proof of Theorem 5.3. □

REMARK 5.2. Recall that in the parametric situation, the risk $\mathbf{E}_{\theta^*} |N_{k^*} \times \mathscr{K}(\widetilde{\theta}^{(k^*)}, \theta^*)|^{1/2}$ of $\widetilde{\theta}^{(k^*)}$ is bounded by $\tau_{1/2}$ (cf. Theorem 2.2). In the non-parametric situation, the result is only slightly worse: the value $\tau_{1/2}$ is replaced by $\sqrt{\tau_1 e^\Delta}$, which takes into account the modeling bias. There is also an additional term proportional to $\sqrt{\overline{\mathfrak{z}}_{k^*}}$, which can be considered as the payment for adaptation. Due to Theorem 5.1, $\overline{\mathfrak{z}}_{k^*}$ is bounded from above by $\mathfrak{z}_K + \iota(K - k^*)$. By Theorem 5.1, $K$ is only logarithmic in the sample size $n$.

Therefore, the risk of the aggregated estimate corresponds to the best possible risk among the family $\{\widetilde{\theta}^{(k)}\}$ for the choice $k = k^*$ up to a logarithmic factor. Lepski, Mammen and Spokoiny [8] established a similar result in the regression setting for the pointwise adaptive Lepski procedure. Combining the result of Corollary 5.5 with Theorem 2.7 yields the rate of adaptive estimation $(n^{-1} \log n)^{1/(2+d)}$ under Lipschitz smoothness of the function $f$ and the usual design regularity; see Polzehl and Spokoiny [12] for more details. It was shown that in the problem of pointwise adaptive estimation, this rate is optimal and cannot be improved by any estimation method. This gives an indirect proof of the optimality of our procedure: the factor $\mathfrak{z}_{k^*}$ in the accuracy of estimation cannot be removed or reduced in the rate because otherwise a similar improvement would appear in the rate of estimation.

## APPENDIX: PROOF OF THEOREM 5.1

The proof utilizes the following simple "metric-like" property of $\mathscr{K}^{1/2}(\cdot, \cdot)$.

LEMMA A.1 (Polzehl and Spokoiny [12], Lemma 5.2). *Under condition* (A2), *it holds for every sequence* $\theta_0, \theta_1, \ldots, \theta_m$ *that*

$$\mathscr{K}^{1/2}(\theta_1, \theta_2) \leq \varkappa \{\mathscr{K}^{1/2}(\theta_1, \theta_0) + \mathscr{K}^{1/2}(\theta_2, \theta_0)\},$$

$$\mathscr{K}^{1/2}(\theta_0, \theta_m) \leq \varkappa \{\mathscr{K}^{1/2}(\theta_0, \theta_1) + \cdots + \mathscr{K}^{1/2}(\theta_{m-1}, \theta_m)\}.$$

With the given constants $\mathfrak{z}_k$, define for $k > 1$ the random sets

$$\mathscr{A}_k = \{N_k \mathscr{K}(\widetilde{\theta}^{(k)}, \widetilde{\theta}^{(k-1)}) \leq b \mathfrak{z}_k\}, \qquad \mathscr{A}^{(k)} = \mathscr{A}_2 \cap \cdots \cap \mathscr{A}_k,$$

where $b$ enters into the construction of $K_{\mathrm{ag}}$: $K_{\mathrm{ag}}(t) = 1$ for $t \leq b$.

First note that $\widehat{\theta}^{(k)} = \widetilde{\theta}^{(k)}$ on $\mathscr{A}^{(k)}$ for all $k \leq K$. This fact can be proved by induction on $k$. For $k = 1$ the assertion is trivial because $\widehat{\theta}^{(1)} = \widetilde{\theta}^{(1)}$. Now suppose that $\widehat{\theta}^{(k-1)} = \widetilde{\theta}^{(k-1)}$. It then holds on $\mathscr{A}_k$ that $\mathbf{m}^{(k)} = N_k \mathscr{K}(\widetilde{\theta}^{(k)}, \widehat{\theta}^{(k-1)}) = N_k \mathscr{K}(\widetilde{\theta}^{(k)}, \widetilde{\theta}^{(k-1)}) \leq b \mathfrak{z}_k$ and thus $\gamma_k = K_{\mathrm{ag}}(\mathbf{m}^{(k)}/\mathfrak{z}_k) \geq K_{\mathrm{ag}}(b) = 1$, yielding $\widehat{\theta}^{(k)} = \widetilde{\theta}^{(k)}$.

Therefore, it remains to bound the risk of $\widehat{\theta}^{(k)}$ on the complement $\overline{\mathscr{A}}^{(k)}$ of $\mathscr{A}^{(k)}$. Define $\mathscr{B}_k = \mathscr{A}^{(k-1)} \setminus \mathscr{A}^{(k)}$. On the event $\mathscr{B}_k$, the index $k$ is the first one for which the condition $N_k \mathscr{K}(\widetilde{\theta}^{(k)}, \widetilde{\theta}^{(k-1)}) \leq b \mathfrak{z}_k$ is violated. It is obvious that $\overline{\mathscr{A}}^{(k)} = \bigcup_{l < k} \mathscr{B}_l$. First, we bound the probability $\mathbf{P}_{\theta^*}(\mathscr{B}_l)$. Applying assumption (A3) and Lemma A.1 yields, for every $l$,

$$N_l \mathscr{K}(\widetilde{\theta}^{(l)}, \widetilde{\theta}^{(l-1)}) \leq 2\varkappa^2 N_l \{\mathscr{K}(\widetilde{\theta}^{(l)}, \theta^*) + \mathscr{K}(\widetilde{\theta}^{(l-1)}, \theta^*)\}$$
$$\leq 2\varkappa^2 \{N_l \mathscr{K}(\widetilde{\theta}^{(l)}, \theta^*) + u_0^{-1} N_{l-1} \mathscr{K}(\widetilde{\theta}^{(l-1)}, \theta^*)\}.$$

Therefore, by Theorem 2.2,

$$\mathbf{P}_{\theta^*}(\mathscr{B}_l) \leq \mathbf{P}_{\theta^*}(N_l \mathscr{K}(\widetilde{\theta}^{(l)}, \widetilde{\theta}^{(l-1)}) > b \mathfrak{z}_l) \leq 2 \exp\left(-\frac{u_0 b}{4\varkappa^2} \mathfrak{z}_l\right).$$

On the set $\mathscr{B}_l$, we have $\widehat{\theta}^{(l-1)} = \widetilde{\theta}^{(l-1)}$ and thus for every $k > l$, the aggregated estimate $\widehat{\theta}^{(k)}$ is, by construction, a convex combination of $\widetilde{\theta}^{(l-1)}, \ldots, \widetilde{\theta}^{(k)}$.

Convexity of the Kullback–Leibler divergence with respect to the second argument, the definition of $\widehat{\theta}^{(k)}$ and Lemma A.1 ensure that

$$\mathscr{K}^{1/2}(\widetilde{\theta}^{(k)}, \widehat{\theta}^{(k)})\mathbf{1}(\mathscr{B}_l) \leq \max_{l'=l-1,\ldots,k-1} \mathscr{K}^{1/2}(\widetilde{\theta}^{(k)}, \widetilde{\theta}^{(l')})$$

$$\leq \varkappa \max_{l'=l-1,\ldots,k-1}\{\mathscr{K}^{1/2}(\widetilde{\theta}^{(k)}, \theta^*) + \mathscr{K}^{1/2}(\widetilde{\theta}^{(l')}, \theta^*)\}$$

$$\leq 2\varkappa \max_{l'=l-1,\ldots,k} \mathscr{K}^{1/2}(\widetilde{\theta}^{(l')}, \theta^*).$$

This and Theorem 2.4 imply that for every $r$,

$$\mathbf{E}_{\theta^*}\mathscr{K}^r(\widetilde{\theta}^{(k)}, \widehat{\theta}^{(k)})\mathbf{1}(\mathscr{B}_l) \leq (2\varkappa)^{2r}\mathbf{E}_{\theta^*}\sum_{l'=l-1}^{k} \mathscr{K}^r(\widetilde{\theta}^{(l')}, \theta^*)\mathbf{1}(\mathscr{B}_l)$$

$$\leq (2\varkappa)^{2r}\sum_{l'=l-1}^{k}\mathbf{E}_{\theta^*}^{1/2}\mathscr{K}^{2r}(\widetilde{\theta}^{(l')}, \theta^*)\mathbf{P}_{\theta^*}^{1/2}(\mathscr{B}_l)$$

$$\leq (2\varkappa)^{2r}\tau_{2r}^{1/2}\sum_{l'=l-1}^{k}N_{l'}^{-r}2\exp\Big(-\frac{u_0 b}{8\varkappa^2}\mathfrak{z}l\Big)$$

$$\leq C_1 N_l^{-r}\tau_{2r}^{1/2}\exp(-c_2\mathfrak{z}l)$$

for some fixed constants $C_1$ and $c_2$. Therefore,

$$\mathbf{E}_{\theta^*}\mathscr{K}^r(\widetilde{\theta}^{(k)}, \widehat{\theta}^{(k)}) \leq \sum_{l=2}^{k}\mathbf{E}_{\theta^*}\mathscr{K}^r(\widetilde{\theta}^{(k)}, \widehat{\theta}^{(k)})\mathbf{1}(\mathscr{B}_l) \leq \sum_{l=2}^{k}C_1 N_l^{-r}\tau_{2r}^{1/2}\exp(-c_2\mathfrak{z}l).$$

It remains to check that the choice $\mathfrak{z}_k = a_0 + a_1\log\alpha^{-1} + a_2 r\log(N_K/N_k)$, with properly selected $a_0, a_1$ and $a_2$, provides the required bound $\mathbf{E}_{\theta^*}|N_k\mathscr{K}(\widetilde{\theta}^{(k)}, \widehat{\theta}^{(k)})|^r \leq \alpha\tau_r$.

## REFERENCES

[1] Belomestny, D. and Spokoiny, V. (2006). Spatial aggregation of local likelihood estimates with applications to classification. SFB 649 Discussion Paper 2006-036.

[2] Breiman, L. (1996). Stacked regressions. *Machine Learning* **24** 49–64.

[3] Cai, Z., Fan, J. and Li, R. (2000). Efficient estimation and inference for varying-coefficient models. *J. Amer. Statist. Assoc.* **95** 888–902. MR1804446

[4] Catoni, O. (2004). *Statistical Learning Theory and Stochastic Optimization. Lecture Notes in Math.* **1851**. Springer, Berlin. MR2163920

[5] Fan, J., Farmen, M. and Gijbels, I. (1998). Local maximum likelihood estimation and inference. *J. R. Stat. Soc. Ser. B Stat. Methodol.* **60** 591–608. MR1626013

[6] Fan, J. and Zhang, W. (1999). Statistical estimation in varying coefficient models. *Ann. Statist.* **27** 1491–1518. MR1742497

 [7] JUDITSKY, A. and NEMIROVSKI, A. (2000). Functional aggregation for nonparametric estimation. *Ann. Statist.* **28** 681–712. MR1792783

 [8] LEPSKI, O., MAMMEN, E. and SPOKOINY, V. (1997). Optimal spatial adaptation to inhomogeneous smoothness: An approach based on kernel estimates with variable bandwidth selectors. *Ann. Statist.* **25** 929–947. MR1447734

 [9] LEPSKI, O. and SPOKOINY, V. (1997). Optimal pointwise adaptive methods in nonparametric estimation. *Ann. Statist.* **25** 2512–2546. MR1604408

[10] LI, J. and BARRON, A. (1999). Mixture density estimation. In *Advances in Neural Inforamtion Processing Systems* **12** (S. A. Sola, T. K. Leen and K. R. Müller, eds.). Morgan Kaufmann Publishers, San Mateo, CA.

[11] LOADER, C. R. (1996). Local likelihood density estimation. *Ann. Statist.* **24** 1602–1618. MR1416652

[12] POLZEHL, J. and SPOKOINY, V. (2006). Propagation-separation approach for local likelihood estimation. *Probab. Theory Related Fields* **135** 335–362. MR2240690

[13] RIGOLLET, PH. and TSYBAKOV, A. (2005). Linear and convex aggregation of density estimators. Manuscript.

[14] SPOKOINY, V. (1998). Estimation of a function with discontinuities via local polynomial fit with an adaptive window choice. *Ann. Statist.* **26** 1356–1378. MR1647669

[15] STANISWALIS, J. G. (1989). The kernel estimate of a regression function in likelihood-based models. *J. Amer. Statist. Assoc.* **84** 276–283. MR0999689

[16] TIBSHIRANI, R. and HASTIE, T. J. (1987). Local likelihood estimation. *Amer. Statist. Assoc.* **82** 559–567. MR0898359

[17] TSYBAKOV, A. (2003). Optimal rates of aggregation. In *Computational Learning Theory and Kernel Machines* (B. Scholkopf and M. Warmuth, eds.) 303–313. *Lecture Notes in Artificial Intelligence* **2777**. Springer, Heidelberg.

[18] YANG, Y. (2004). Aggregating regression procedures to improve performance. *Bernoulli* **10** 25–47. MR2044592

WEIERSTRASS INSTITUTE
MOHRENSTRASSE 39
10117 BERLIN
GERMANY
E-MAIL: belomest@wias-berlin.de
        spokoiny@wias-berlin.de